\documentclass{amsart}
\usepackage{amssymb,amscd}
\usepackage[all]{xy}

\theoremstyle{plain}

\theoremstyle{definition}

\theoremstyle{remark}

\newtheorem*{ack}{Acknowledgements}

\newcommand{\Curve}{\mathcal{C}}
\renewcommand{\O}{\mathcal{O}}
\newcommand{\FF}{\mathcal{F}}
\newcommand{\EE}{\mathcal{E}}
\newcommand{\LL}{\mathcal{L}}

\newcommand{\F}{\mathbb{F}}
\newcommand{\Fp}{{\mathbb{F}_p}}
\newcommand{\Fq}{{\mathbb{F}_q}}
\newcommand{\Z}{\mathbb{Z}}
\newcommand{\Q}{\mathbb{Q}}
\newcommand{\R}{\mathbb{R}}
\newcommand{\C}{\mathbb{C}}
\newcommand{\A}{\mathbb{A}}
\renewcommand{\H}{\mathbb{H}}
\renewcommand{\P}{\mathbb{P}}

\newcommand{\m}{\mathfrak{m}}
\newcommand{\n}{\mathfrak{n}}
\renewcommand{\a}{\mathfrak{a}}
\newcommand{\p}{\mathfrak{p}}

\newcommand{\sha}{{\hbox to 10pt{\rlap{\hskip2.8pt\vrule
height6pt\hskip1.6pt\vrule height6pt\hskip1.6pt
\vrule height6pt}\hskip1pt\vrule height0.8pt width 8pt\hskip1pt}}}

\newcommand{\into}{\hookrightarrow}

\newcommand{\tensor}{\otimes}
\newcommand{\compose}{\circ}
\def\nodiv{\mathrel{\mathchoice{\not|}{\not|}{\kern-.2em\not\kern.2em|}{\kern-.2em\not\kern.2em|}}}
\newcommand{\PGL}{\mathrm{PGL}}

\newcommand{\GL}{\mathrm{GL}}

\DeclareMathOperator{\res}{Res}

\DeclareMathOperator{\ord}{ord}
\DeclareMathOperator{\rk}{Rank}

\DeclareMathOperator{\Hom}{Hom}

\DeclareMathOperator{\gal}{Gal}
\DeclareMathOperator{\spec}{Spec}
\DeclareMathOperator{\en}{End}
\DeclareMathOperator{\mor}{Mor}

\theoremstyle{plain}
\newtheorem{conj}[subsection]{Conjecture}
\newtheorem{thm}[subsubsection]{Theorem}

\numberwithin{equation}{section}

\newcommand{\invlim}{\underleftarrow{\lim}}

\newcommand{\Fpbar}{{\overline{\mathbb{F}}_p}}
\newcommand{\Fqbar}{{\overline{\mathbb{F}}_q}}

\newcommand{\Fqq}{{\mathbb{F}}_{q^2}}
\newcommand{\Fr}{{\mathbb{F}}_{r}}

\newcommand{\Qbar}{{\overline{\mathbb{Q}}}}
\newcommand{\Ql}{{\mathbb{Q}_\ell}}
\newcommand{\Qlbar}{{\overline{\mathbb{Q}}_\ell}}

\newcommand{\ratto}{{\dashrightarrow}}
\newtheorem*{thm*}{Theorem}

\def\hyp#1{{\advance\hsize by -1in\parindent=0pt
    \vtop{#1}}}

\begin{document}
\title[Elliptic curves and analogies]{Elliptic
curves and analogies between number fields and function fields}
\author{Douglas Ulmer}
\address{Department of Mathematics \\ 
         University of Arizona \\ Tucson,
         AZ  85721}
\email{ulmer@math.arizona.edu}
\thanks{This paper is based upon work supported by the National
Science Foundation under Grant No. DMS 0070839}
\date{June 2, 2003}
\begin{abstract}
  The well-known analogies between number fields and function fields
  have led to the transposition of many problems from one domain to
  the other.  In this paper, we will discuss traffic of this sort, in
  both directions, in the theory of elliptic curves.  In the first
  part of the paper, we consider various works on Heegner points and
  Gross-Zagier formulas in the function field context; these works
  lead to a complete proof of the conjecture of Birch and
  Swinnerton-Dyer for elliptic curves of analytic rank at most 1 over
  function fields of characteristic $>3$.  In the second part of the
  paper, we will review the fact that the rank conjecture for elliptic
  curves over function fields is now known to be true, and that the
  curves which prove this have asymptotically maximal rank for their
  conductors.  The fact that these curves meet rank bounds suggests a
  number of interesting problems on elliptic curves over number
  fields, cyclotomic fields, and function fields over number fields.
  These problems are discussed in the last four sections of the paper.
\end{abstract}
\maketitle

\section{Introduction}

The purpose of this paper is to discuss some work on elliptic curves
over function fields inspired by the Gross-Zagier theorem and some new
ideas about ranks of elliptic curves from the function field case
which I hope will inspire work over number fields.

We begin in Section~\ref{s:BSD} by reviewing the statement of and
current state of knowledge on the conjecture of Birch and
Swinnerton-Dyer for elliptic curves over function fields.  Then in
Section~\ref{s:ffGZ} we discuss various works by R\"uck and Tipp,
P\'al, and Longhi on function field analogues of the Gross-Zagier
formula and some related work by Brown.  We also explain how suitably
general Gross-Zagier formulas together with my
``geometric non-vanishing'' results lead to a theorem
of the form: the Birch and Swinnerton-Dyer conjecture for elliptic
curves over function fields of curves over finite fields of
characteristic $>3$ holds for elliptic curves with analytic rank at
most 1.

In Sections~\ref{s:ffRanks} and \ref{s:bounds} we move beyond rank one
and explain that the rank conjecture holds for elliptic curves over
function fields: there are (non-isotrivial) elliptic curves with
Mordell-Weil group of arbitrarily large rank.  Moreover, these curves
meet an asymptotic bound due to Brumer for the rank in terms of the
conductor.  So in the function field case, we know precisely the
asymptotic growth of ranks of elliptic curves ordered by the size of
their conductors.  In fact, there are two bounds, one arithmetic, the
other geometric, and both are sharp.

The rest of the paper is devoted to explaining some interesting problems
suggested by the existence and sharpness of these two types of rank
bounds.  In Section~\ref{s:nfRanks} we make a conjecture which says
roughly that Mestre's bound on the ranks of elliptic curves over $\Q$
and suitable generalizations of it over number fields are asymptotically
sharp.  Next, we note that the Mestre bound and even more so the Brumer
bound are (or rather can be reformulated as) algebraic statements.
For example, the Brumer bound can be interpreted as a statement about
the eigenvalues of Frobenius on \'etale cohomology.  It is therefore
natural to ask for an algebraic proof; reformulating the bounds into
statements that might admit an algebraic proof leads to some
interesting questions which are explained in
Section~\ref{s:algRankBounds}.

Finally, in Sections~\ref{s:cycRanks} and~\ref{s:higherdffRanks} we
discuss possible rank bounds over cyclotomic fields and over function
fields over number fields.  More precisely, we discuss pairs of ranks
bounds, one ``arithmetic'' the other ``geometric,'' for pairs of
fields like $\Q^{p-\rm cyc}/\Q$ or $\Qbar(\Curve)/\Q(\Curve)$ where
$\Curve$ is a curve over $\Q$.  In both cases, one rank bound is known
(arithmetic in the first case, geometric in the second) and the other
bound has yet to be considered.

\begin{ack}
It is a pleasure to thank Gautam Chinta, Henri Darmon, Mihran
Papikian, Joe Silverman, Dinesh Thakur, Adrian Vasiu and the referee
for their comments, corrections, and references to the literature.
\end{ack}

\section{Review of the Birch and Swinnerton-Dyer conjecture over
  function fields}\label{s:BSD} 
We assume that the reader is familiar with elliptic curves over number
fields, but perhaps not over function fields, and so in this
preliminary section we set up some background and review the Birch and
Swinnerton-Dyer conjecture.  For many more details, examples, etc., we
refer to \cite{UlmerS}.

Let $\Curve$ be a smooth, geometrically connected, projective curve over
a finite field $\Fq$ and set $F=\Fq(\Curve)$.  Let $E$ be an elliptic
curve over $F$, i.e., a curve of genus one defined as usual by an
affine Weierstrass equation
\begin{equation*}
y^2+a_1xy+a_3y=x^3+a_2x^2+a_4x+a_6
\end{equation*}
($a_i\in F$) with the point at infinity $[0,1,0]$ as origin; the
discriminant $\Delta$ and $j$-invariant are given by the usual
formulas (see, e.g., \cite{TateAlg}) and we of course assume that
$\Delta\neq0$.  We say that $E$ is {\it constant\/} if it is 
defined over $\Fq$, i.e., if it is possible to choose the Weierstrass
model so that the $a_i\in\Fq$.  Equivalently, $E$ is constant if
there exists an elliptic curve $E_0$ defined over $\Fq$ such that
$E\cong E_0\times_{\spec\Fq}\spec F$.  In this case we say that $E$ is
based on $E_0$.  We say that $E$ is {\it isotrivial\/} if it becomes
isomorphic to a constant curve after a finite extension of $F$; this
is easily seen to be equivalent to the condition $j(E)\in\Fq$.  Finally,
we say that $E$ is {\it non-isotrivial\/} if $j(E)\not\in\Fq$.

Let $\n$ be the conductor of $E$.  This is an effective divisor on
$\Curve$ which is divisible only by the places where $E$ has bad
reduction.  More precisely, $v$ divides $\n$ to order 1 at places
where $E$ has multiplicative reduction and to order at least 2 at
places where $E$ has additive reduction and to order exactly 2 at
these places if the characteristic of $F$ is $>3$. The reduction,
exponent of conductor, and minimal model of $E$ at places of $F$ can
be computed by Tate's algorithm \cite{TateAlg}.

The Mordell-Weil theorem holds for $E$, namely $E(F)$ is a finitely
generated abelian group.  This can be proven in a manner entirely
analogous to the proof over number fields, using Selmer groups and
heights, or by more geometric methods; see \cite{Neron}.  Also, both
the rank conjecture (that for a fixed $F$, the rank of $E(F)$ can be
arbitrarily large) and the torsion conjecture (that there is a bound
on the order of the torsion subgroup of $E(F)$ depending only on the
genus of $F$) are known to be true in this context.  For the rank
conjecture, see \cite{UlmerR} and Section~\ref{s:ffRanks} below.  The
torsion conjecture was proven by Levin \cite{Levin}, who showed that
there is an explicit bound of the form $O(\sqrt{g}+1)$ for the order of
the torsion subgroup of a non-isotrivial elliptic curve over $F$,
where $g$ is the genus of $F$.  More recently, Thakur \cite{Thakur}
proved a variant bounding the order of torsion in terms of the
``gonality'' of $\Curve$, i.e., the smallest degree of a non-constant
map to $\P^1$.

The $L$-function of $E$ is defined by the Euler product
\begin{align*}
L(E/F,s)=&\prod_{v\nodiv\n}\left(1-a_vq_v^{-s}+q_v^{1-2s}\right)^{-1}\cr
&\times\prod_{v|\n}\begin{cases}
(1-q_v^{-s})^{-1}&
     \text{if $E$ has split multiplicative reduction at $v$}\cr
(1+q_v^{-s})^{-1}&
     \text{if $E$ has non-split multiplicative reduction at $v$}\cr
1&\text{if $E$ has additive reduction at $v$.}\cr
               \end{cases}
\end{align*}
Here $q_v$ is the cardinality of the residue field $\F_v$ at $v$ and
the number of points on the reduced curve is $\#E(\F_v)=q_v+1-a_v$.
The product converges absolutely in the half-plane $\Re s>3/2$,
has a meromorphic continuation to the $s$ plane, and
satisfies a functional equation for $s\mapsto 2-s$.  If $E$ is not
constant, then $L(E/F,s)$ is a polynomial in $s$ of degree
$4g-4+\deg\n$ and thus an entire function of $s$.  (All this comes from
Grothendieck's analysis of $L$-functions.  See the last section of
\cite{MilneEC} for more details.)

The Birch and Swinnerton-Dyer conjecture in this context
asserts that
\begin{equation*}
\rk E(F)\mathrel{\mathop=^?}\ord_{s=1}L(E/F,s)
\end{equation*}
and, setting $r=\ord_{s=1}L(E/F,s)$, that the leading coefficient is
\begin{equation*}
\frac1{r!}L^{(r)}(E/F,1)\mathrel{\mathop=^?}
\frac{|\sha| R\tau}{|E(F)_{tor}|^2}
\end{equation*}
where $\sha$ is the Tate-Shafarevitch group, $R$ is a regulator
constructed from heights of a set of generators of $E(F)$, and $\tau$
is a certain Tamagawa number (an analogue of a period).  We will not
enter into the details of the definitions of these objects since they
will play little role in what follows; see \cite{TateBSD} for more details.

Much more is known about this conjecture in the function field case
than in the number field case.  Indeed, we have
\begin{equation}\label{eq:rk<ord}
\rk E(F)\le\ord_{s=1}L(E/F,s)
\end{equation}
and the following assertions are equivalent:
\begin{enumerate}
\item Equality holds in \ref{eq:rk<ord}.
\item The $\ell$ primary part of $\sha$ is finite for any one prime
$\ell$ ($\ell=p$ is allowed).
\item $\sha$ is finite.
\end{enumerate}
Moreover, if these equivalent conditions are satisfied, then the
refined conjecture on the leading coefficient of the $L$-series is
true.  The ``prime-to-$p$'' part of this was proven by Artin and Tate
\cite{TateBSD}.  More precisely, they showed that equality holds in
\ref{eq:rk<ord} if and only if the $\ell$ primary part of $\sha$ is
finite for any one prime $\ell\neq p$ if and only if the $\ell$
primary part of $\sha$ is finite for every $\ell\neq p$, and that if these
conditions hold, the refined formula is correct up to a power of $p$.
Milne proved the stronger statement above in \cite{MilneBSD} for
$p\neq2$; due to later improvements in $p$-adic cohomology, his
argument applies essentially verbatim to the case $p=2$ as well.

These results were obtained by considering the elliptic surface
$\EE\to\Curve$ attached to $E$, which can be characterized as the
unique smooth, proper surface over $\Fq$ admitting a flat and
relatively minimal morphism to $\Curve$, with generic fiber $E/F$.
Another key ingredient is Grothendieck's analysis of $L$-functions,
which gives a cohomological interpretation of the $\zeta$-function of
$\EE$ and the $L$-function of $E$.

Equality in \ref{eq:rk<ord}, and therefore the full Birch and
Swinnerton-Dyer conjecture, is known to hold in several cases (but
certainly not the general case!): If it holds for $E/K$ where $K$ is a
finite extension of $F$, then it holds for $E/F$ (this is elementary);
it holds for constant, and thus isotrivial, $E$ (this follow from
\cite{TateEnd}); and it holds for several cases most easily described
in terms of $\EE$, namely if $\EE$ is a rational surface (elementary),
a $K3$ surface \cite{ArtinSD}, or if $\EE$ is dominated by a product
of curves (see \cite{TateM}).  The rational and $K3$ cases are
essentially those where the base field $F$ is $\Fq(t)$ and the
coefficients $a_i$ in the defining Weierstrass equation of $E$ have
small degree in $t$.

\numberwithin{equation}{subsection}
\section{Function field analogues of the Gross-Zagier theorem}\label{s:ffGZ}
In this section we will give some background on modularity and Heegner
points and then discuss various works on Gross-Zagier formulas in the
function field context.  Our treatment will be very sketchy, just
giving the main lines of the arguments, but we will give precise
references where the reader may find the complete story.
Throughout, we fix a smooth, proper, geometrically connected curve
$\Curve$ over a finite field $\Fq$ of characteristic $p$ and we set
$F=\Fq(\Curve)$.

\subsection{Two versions of modularity}
Recall that for elliptic curves over $\Q$ there are two
(not at all trivially!) equivalent statements expressing the property that
an elliptic curve $E$ of conductor $N$ is modular:
\begin{enumerate}
\item There exists a modular form $f$ (holomorphic of weight 2 and
level $\Gamma_0(N)$) such that $L(E,\chi,s)=L(f,\chi,s)$ for all Dirichlet
characters $\chi$.
\item There exists a non-constant morphism 
$X_0(N)\to E$ (defined over $\Q$).
\end{enumerate}
(We note that over $\Q$, the equalities $L(E,\chi,s)=L(f,\chi,s)$ in
(1) are implied by the {\it a priori\/} weaker statement that
$L(E,s)=L(f,s)$.  But this implication fails over higher degree number
fields and over function fields and it is the stronger assertion that
we need.)

In the next two subsections we will explain the analogues of these two
statements in the function field context.  In this case, the
relationship between the two statements is a little more complicated
than in the classical case.  For example, the relevant automorphic
forms are complex valued and thus are not functions or sections of
line bundles on the analogue of $X_0(\n)$, which is a curve over $F$.
Nevertheless, analogues of both modularity statements are theorems in
the function field case.

\subsection{Analytic modularity}\label{ss:AnalyticMod}
We begin with (1).  Let $\A_F$ be the ad\`ele ring of $F$ and
$\O_F\subset\A_F$ the subring of everywhere integral ad\`eles.  Then
for us, automorphic forms on $\GL_2$ over $F$ are functions on
$\GL_2(\A_F)$ which are invariant under left translations by
$\GL_2(F)$ and under right translations by a finite index subgroup $K$
of $\GL_2(\O_F)$.  In other words, they are functions on the double 
coset space
\begin{equation}\label{eqn:doublecoset}
\GL_2(F)\backslash \GL_2(\A_F)/K.
\end{equation}
These functions may take values in any field of
characteristic zero; to fix ideas, we take them with values in $\Qbar$
and we fix embeddings of $\Qbar$ into $\C$ and into $\Qlbar$ for some
$\ell\neq p$.  The subgroup $K$ is the analogue of the level in the
classical setting and the most interesting case is when $K$ is one of
the analogues $\Gamma_0(\m)$ or $\Gamma_1(\m)$ of the Hecke congruence
subgroups where $\m$ is an effective divisor on $\Curve$.  If
$\psi:\A^\times/F^\times\to\Qlbar^\times$ is an id\`ele class
character and $f$ is an automorphic form, we say $f$ has central
character $\psi$ if $f(zg)=\psi(z)f(g)$ for all $z\in
Z(\GL_2(\A_F))\cong\A_F^\times$ and all $g\in\GL_2(\A_F)$.  The
central character plays the role of weight: When $k$ is a positive
integer and $\psi(z)=|z|^{-k}$ (where $|\cdot|$ is the ad\`elic norm),
$f$ is analogous to a classical modular form of weight $k$.  The basic
reference for this point of view is \cite{WeilDS}; see Chapter III for
definitions and first properties.  For a more representation-theoretic
point of view, see \cite{JL}.

If we single out a place $\infty$ of $F$ and assume that
$K=\Gamma_0(\infty\n)$ where $\n$ is prime to $\infty$, then there is
an analogue of the description of classical modular forms as functions
on the upper half plane.  Namely, an automorphic form $f$ may be
viewed as a function (or a section of line bundle if the $\infty$
component of $\psi$ is non-trivial) on a finite number of copies of
the homogeneous space $\PGL_2(F_\infty)/\Gamma_0(\infty)$ which has
the structure of an oriented tree.  (Compare with
$\PGL_2(\R)/O_2(\R)\cong\H$.)  The corresponding functions are
invariant under certain finite index subgroups of
$\GL_2(A)\subset\GL_2(F_\infty)$ where $A\subset F$ is the ring of
functions regular outside $\infty$.  The various copies of the tree
are indexed by a generalized ideal class group of $A$.  This point of
view is most natural when $F=\Fq(t)$ and $\infty$ is the standard
place $t=\infty$, in which case there is just one copy of the tree and
this description is fairly canonical.  In the general case, there are
several copies of the tree and choices must be made to identify
automorphic forms with functions on trees.  Using this description (or
suitable Hecke operators) one may define the notion of a form being
``harmonic'' or ``special'' at $\infty$.  Namely, sum of the values
over edges with a fixed terminus should be zero.  This is an analogue
of being holomorphic.  See \cite[Chap.~5]{DH}, \cite[Chap.~4]{GekRev}, or
\cite[Chap.~2]{vanderPutReversat} for details.

Automorphic forms have Fourier expansions, with coefficients naturally
indexed by effective divisors on $\Curve$.  There are Hecke operators,
also indexed by effective divisors on $\Curve$ and the usual
connection between eigenvalues of Hecke operators and Fourier
coefficients of eigenforms holds.  There is a notion of cusp form and
for a fixed $K$ and $\psi$ the space of cusp forms is finite
dimensional.  An automorphic form $f$ gives rise to an $L$-function
$L(f,s)$, which is a complex valued function of a complex variable
$s$. If $f$ is a cuspidal eigenform, this $L$-function has an Euler
product, an analytic continuation to an entire function of $s$, and
satisfies a functional equation.  See \cite{WeilDS} for all of this
except the finite dimensionality, which follows easily from reduction
theory.  See \cite[Chap.~II]{Serre} for the finite dimensionality when
$F=\Fq(t)$ and \cite{HLW} for an explicit dimension formula in the
general case.

The main theorem of \cite{WeilDS} is a ``converse'' theorem which says
roughly that a Dirichlet series with a suitable analytic properties is
the $L$-function of an automorphic form on $\GL_2$.  (The function
field case is Theorem 3 of Chapter VII.) The most important
requirement is that sufficiently many of the twists of the given
Dirichlet series by finite order characters should satisfy functional
equations.  This result was also obtained by representation theoretic
methods in \cite{JL}.  Also, see \cite{Li} for an improved version,
along the lines of \cite{WeilDS}.

Now let $E$ be an elliptic curve over $F$.  By Grothendieck's analysis
of $L$-functions, we know that the Dirichlet series $L(E,s)$ is
meromorphic (entire if $E$ is non-isotrivial) and its twists satisfy
functional equations.  In \cite[9.5-9.7]{DeligneC}, Deligne verified
the hypotheses of Weil's converse theorem.  The main point is to check
that the functional equations given by Grothendieck's theory are the
same as those required by Weil.  The form $f_E$ associated to $E$ is
characterized by the equalities $L(E,\chi,s)=L(f_E,\chi,s)$ for all
finite order id\`ele class characters $\chi$.  It is an eigenform for
the Hecke operators and is a cusp form if $E$ is non-isotrivial.  Its
level is $\Gamma_0(\m)$ where $\m$ is the conductor of $E$ and it has
central character $|\cdot|^{-2}$ (i.e., is analogous to a form of
weight 2).  If $E$ has split multiplication at $\infty$, then $f_E$ is
special at $\infty$.  The construction of $f_E$ from $E$ is the
function field analogue of (1) above.

\subsection{Geometric modularity}
We now turn to Drinfeld modules and (2).  There is a vast literature
on Drinfeld modules and we will barely scratch the surface.  The
primary reference is \cite{Drinfeld} and there are valuable surveys
in \cite{DH} and \cite{AB}.  

Fix a place $\infty$ of $F$ and define $A$ to be the ring of elements
of $F$ regular away from $\infty$.  Let $F_\infty$ denote the
completion of $F$ at $\infty$ and $C$ the completion of the algebraic
closure of $F_\infty$.  The standard example is when $F=\Fq(t)$,
$\infty$ is the standard place $t=\infty$, and $A=\Fq[t]$.

Let $k$ be a ring of characteristic $p$ equipped with a homomorphism
$A\to k$.  Let $k\{\tau\}$ be the ring of non-commutative polynomials
in $\tau$, with commutation relation $\tau a=a^p\tau$.  There is a
natural inclusion $\epsilon:k\into k\{\tau\}$ with left inverse
$D:k\{\tau\}\to k$ defined by $D(\sum_na_n\tau^n)=a_0$.  If $R$ is any
$k$-algebra, we may make the additive group of $R$ into a module over
$k\{\tau\}$ by defining $(\sum_na_n\tau^n)(x)=\sum_na_nx^{p^n}$.

A Drinfeld module over $k$ (or elliptic module as Drinfeld called
them) is a ring homomorphism $\phi:A\to k\{\tau\}$ whose image is not
contained in $k$ and such that $D\compose\phi:A\to k$ is the given
homomorphism.  The characteristic of $\phi$ is by definition the
kernel of the homomorphism $A\to k$, which is a prime ideal of $A$.
It is convenient to denote the image of $a\in A$ by $\phi_a$ rather
than $\phi(a)$.  If $A=\Fq[t]$ then $\phi$ is determined by $\phi_t$,
which can be any element of $k\{\tau\}$ of positive degree with
constant term equal to the image of $t$ under $A\to k$.  For a general
$A$ and $k$ equipped with $A\to k$ there may not exist any Drinfeld
modules and if they do exist, they may not be easy to find.
As above, a Drinfeld module $\phi$ turns any $k$-algebra into an
$A$-module by the rule $a\cdot x=\phi_a(x)$.

It turns out that $a\mapsto\phi_a$ is always injective and there
exists a positive integer $r$, the rank of $\phi$, such that
$p^{\deg_\tau(\phi_a)}=|a|_\infty^r=\#(A/a)^r$. If $\phi$ and $\phi'$
are Drinfeld modules over $k$, a homomorphism $u:\phi\to\phi'$ is by
definition an element $u\in k\{\tau\}$ such that $u\phi_a=\phi'_au$
for all $a\in A$ and an isogeny is a non-zero homomorphism.  Isogenous
Drinfeld modules have the same rank and characteristic.  See
\cite[\S2]{Drinfeld} or \cite[Chap.~1]{DH}.

We will only consider Drinfeld modules of rank 2.  These objects are
in many ways analogous to elliptic curves.  For example, if $k$ is an
algebraically closed field and $\p\subset A$ is a prime ideal, then we have an isomorphism
of $A$-modules
\begin{equation*}
\phi[\p](k):=\{x\in k|\phi_a(x)=0\text{ for all }a\in\p\}\cong(A/\p)^e
\end{equation*}
where $0\le e\le2$ and $e=2$ if the characteristic of $\phi$ is
relatively prime to $\p$.  A second analogy occurs with endomorphism
rings: $\en(\phi)$, the ring of endomorphisms of $\phi$, is isomorphic
as $A$-module to either $A$, an $A$-order in an ``imaginary''
quadratic extension $K$ of $F$, or an $A$-order in a quaternion
algebra over $F$.  Here ``imaginary'' means that the place $\infty$ of
$F$ does not split in $K$ and an $A$-order in a division algebra $D$
over $F$ is an $A$-subalgebra $R$ which is projective of rank 4 as
$A$-module.  The quaternion case can occur only if the characteristic
of $\phi$ is non-zero, in which case the quaternion algebra is
ramified precisely at $\infty$ and the characteristic of $\phi$.  A
third analogy is the analytic description of Drinfeld modules over
$C$: giving a Drinfeld module of rank 2 over $C$ up to isomorphism is
equivalent to giving a rank 2 $A$-lattice in $C$ up to homothety by
elements of $C^\times$.  If $\phi$ corresponds to the lattice
$\Lambda$, there is a commutative diagram
\begin{equation*}
\xymatrix{0\ar[r]&\Lambda\ar[r]\ar[d]^a&C\ar[r]^{\exp_\Lambda}\ar[d]^a
&C\ar[d]^{\phi_a}\ar[r]&0\\ 
0\ar[r]&\Lambda\ar[r]&C\ar[r]^{\exp_\Lambda}&C\ar[r]&0}
\end{equation*}
where $\exp_\Lambda:C\to C$ is the Drinfeld exponential associated to
$\Lambda$.  See \cite[\S\S2-3]{Drinfeld} or \cite[Chaps. 1-2]{DH}.

There is a natural generalization of all of the above to Drinfeld
modules over schemes of characteristic $p$.  Given an effective
divisor $\n$ on $\Curve$ relatively prime to $\infty$ (or
equivalently, a non-zero ideal of $A$), there is a notion of ``level
$\n$ structure'' on a Drinfeld module.  Using this notion, one may
construct a moduli space $Y_0(\n)$ (a scheme if $\n$ is non-trivial
and a stack if $\n$ is trivial) parameterizing Drinfeld modules of
rank 2 with level $\n$ structure, or equivalently, pairs of rank 2
Drinfeld modules connected by a ``cyclic $\n$-isogeny''
$u:\phi\to\phi'$.  (The notions of level $\n$ structure and cyclic
isogeny are somewhat subtle and a significant advance over the naive
notions.  Analogues of Drinfeld's notions were used in
\cite{KatzMazur} to completely analyze the reduction of classical
modular curves at primes dividing the level.)  The curve $Y_0(\n)$ is
smooth and affine over $F$ and may be completed to a smooth, proper
curve $X_0(\n)$.  The added points (``cusps'') can be interpreted in
terms of certain degenerations of Drinfeld modules.  The curve
$X_0(\n)$ carries many of the structures familiar from the classical
case, such as Hecke correspondences (indexed by effective divisors on
$\Curve$) and Atkin-Lehner involutions.  See \cite[\S5]{Drinfeld} and
\cite[Chap.~1, \S6]{DH}.  The construction of the moduli space (or
stack) is done very carefully in \cite[Chap.~1]{Laumon} and the
interpretation of the cusps is given in \cite{vanderPutTop}.

The analytic description of Drinfeld modules over $C$ yields an
analytic description of the $C$ points of $Y_0(\n)$.  Namely, let
$\Omega$ denote the Drinfeld upper half plane:
$\Omega=\P^1(C)\setminus P^1(F_\infty)$.  Then $Y_0(\n)(C)$ is
isomorphic (as rigid analytic space) to a union of quotients of
$\Omega$ by finite index subgroups of $\GL_2(A)$.  The components of
$Y_0(\n)(C)$ are indexed by a generalized ideal class group of $A$.
More adelically, we have an isomorphism 
\begin{equation*}
Y_0(\n)\cong\GL_2(F)\backslash\left(GL_2(\A_F^f)\times\Omega\right)
/\Gamma_0(\n)^f
\end{equation*}
where $\A_F^f$ denotes the ``finite ad\`eles'' of $F$, namely the
ad\`eles with the component at $\infty$ removed, and similarly with
$\Gamma_0(\n)^f$.  See \cite[\S6]{Drinfeld} or \cite[Chap.~3]{DH}. 

This description reveals a close connection between $Y_0(\n)$ and the
description of automorphic forms as functions on trees
(cf.~\ref{eqn:doublecoset}).  Namely, there is a map between
the Drinfeld upper half plane $\Omega$ and a geometric realization of
the tree $\PGL_2(F_\infty)/\Gamma_0(\infty)$.  Using this, Drinfeld
was able to analyze the \'etale cohomology of $X_0(\n)$ as a module
for $\gal(\overline{F}/F)$ and the Hecke operators, in terms of
automorphic forms of level $\Gamma_0(\n\infty)$ which are special at
$\infty$.  (Drinfeld used an {\it ad hoc\/} definition of \'etale
cohomology; for a more modern treatment, see
\cite{vanderPutReversat}.)  This leads to one form of the Drinfeld
reciprocity theorem: if $f$ is an eigenform of level
$\Gamma_0(\n\infty)$ which is special at $\infty$, then there exists a
factor $A_f$ of the Jacobian $J_0(\n)$ of $X_0(\n)$, well-defined up
to isogeny, such that
\begin{equation*}
L(A_f,\chi,s)=\prod_{\sigma:\Q(f)\into\C}L(f^\sigma,\chi,s)
\end{equation*}
for all finite order id\`ele class characters $\chi$ of $F$.  Here the
product is over all embeddings of the number field generated by the
Fourier coefficients of $f$ into $\C$.  If the Hecke eigenvalues of
$f$ are rational integers, then $A_f$ is an elliptic curve, and if $f$
is a new, then $E$ has conductor $\n\infty$ and is split
multiplicative at $\infty$.  See \cite[\S\S 10-11]{Drinfeld},
\cite[Chaps.~4-5]{DH}, and \cite[Chap.~8]{GekRev}.

So, starting with an elliptic curve $E$ over $F$ of level
$\m=\n\infty$ which is split multiplicative at $\infty$, Deligne's
theorem gives us an automorphic form $f_E$ on $\GL_2$ over $F$ of
level $\m$ which is special at $\infty$ and which has integer Hecke
eigenvalues.  From $f_E$, Drinfeld's construction gives us an isogeny
class of elliptic curves $A_{f_E}$ appearing in the Jacobian of
$X_0(\n)$.  Moreover, we have equalities of $L$-functions:
\begin{equation*}
L(E,\chi,s)=L(f_E,\chi,s)=L(A_{f_E},\chi,s).
\end{equation*}
But Zarhin proved that the $L$-function of an abelian
variety $A$ over a function field (by which we mean the collection of all
twists $L(A,\chi,s)$) determines its isogeny class.  See
\cite{Zarhin} for the case $p>2$ and \cite[XXI.2]{MoretBailly} for a
different proof that works in all characteristics.  This means that
$E$ is in the class $A_{f_E}$ and therefore we have a non-trivial
modular parameterization $X_0(\n)\to E$.

In \cite[Chap.~9]{GekRev} Gekeler and Reversat completed this picture
by giving a beautiful analytic construction of $J_0(N)(C)$ and of the
analytic parameterization $X_0(\n)(C)\to E(C)$.  This is the
analogue of the classical parameterization of an elliptic curve by
modular functions.  Recently, Papikian has studied the degrees of
Drinfeld modular parameterizations and proved the analogue of the
degree conjecture.  See \cite{Papikian} and forthcoming publications.

\subsection{Heegner points and Brown's work}
It was clear to the experts from the beginning that Heegner points,
the Gross-Zagier formula, and Kolyvagin's work could all be extended
to the function field case, using the Drinfeld modular
parameterization, although people were reluctant to do so.  The first
efforts in this direction were made by Brown in \cite{Brown}.

Fix as usual $F$, $\infty$, $A$, and $\n$, so we have the Drinfeld
modular curve $X_0(\n)$.  Let $K/F$ be an imaginary quadratic
extension and let $B$ be an $A$-order in $K$.  A Drinfeld-Heegner
point with order $B$ (or Heegner point for short) is by definition a
point on $X_0(\n)$ corresponding to a pair $\phi\to\phi'$ connected by
a cyclic $\n$ isogeny such that $\en(\phi)=\en(\phi')=B$.  These will
exist if and only if there exists a proper ideal $\n'$ of $B$ (i.e.,
one such that $\{b\in K|b\n'\subset\n'\}=B$) with $B/\n'\cong A/\n$.
The simplest situation is when every prime dividing $\n$ splits in $K$
and $B$ is the maximal $A$-order in $K$, i.e., the integral closure of
$A$ in $K$.  Assuming $B$ has such an ideal $\n'$, we may construct
Heegner points using the analytic description of Drinfeld modules over
$C$ as follows.  If $\a\subset B$ is a non-zero proper ideal then the
pair of Drinfeld modules $\phi$ and $\phi'$ corresponding to the
lattices $\n'\a$ and $\a$ in $K\into C$ satisfy
$\en(\phi)=\en(\phi')=B$ and they are connected by a cyclic
$\n$-isogeny.  The corresponding point turns out to depend only on
$B$, $\n'$, and the class of $\a$ in $\text{Pic}(B)$ and it is defined
over the ring class field extension $K_B/K$ corresponding to $B$ by
class field theory.  The theory of complex multiplication of Drinfeld
modules implies that $\gal(K_B/K)\cong\text{Pic}(B)$ acts on the
Heegner points through its natural action on the class of $\a$.
Applying an Atkin-Lehner involution to a Heegner point is related to
changing the choice of ideal $\n'$ over $\n$.  All of this is
discussed in \cite[\S2]{Brown} in the context where $A=\Fq[t]$.

Taking the trace from $K_B$ to $K$ of a Heegner point and subtracting
a suitable multiple of a cusp, we get a $K$-rational divisor of degree
0, and so a point $J_0(\n)(K)$.  We write $Q_K$ for the point so
constructed when $K$ is an imaginary quadratic extension of $F$ in
which every prime dividing $\n$ splits and $B$ is the maximal $A$-order in
$K$.  The point $Q_K$ is well-defined, independently of the other
choices ($\n'$ and $\a$), up to a torsion point of $J_0(\n)(K)$.  If
$E$ is an elliptic curve over $F$ of level $\n\infty$ with split
multiplicative reduction at $\infty$, then using the modular
parameterization discussed above one obtains a point $P_K\in E(K)$,
well-defined up to torsion.

Brown purports to prove, by methods analogous to those of
Kolyvagin~\cite{Kolyvagin}, that if $P_K$ is non-torsion, then the
Tate-Shafarevitch group of $E$ is finite and the rank of $E(K)$ is
one.  (He gives an explicit annihilator of the $\ell$-primary part of
$\sha$ for infinitely many $\ell$.)  As we have seen, this implies
that the Birch and Swinnerton-Dyer conjecture holds for $E$ over $K$.

Unfortunately, Brown's paper is marred by a number of errors, some
rather glaring.  For example, the statement of the main theorem is not
in fact what is proved and it is easily seen to be false if taken
literally.  Also, he makes the strange hypothesis that $q$, the number
of elements in the finite ground field, is not a square.  The source
of this turns out to be a misunderstanding of quadratic reciprocity in
the proof of his Corollary~3.4. In my opinion, although something like
what Brown claims can be proved by the methods in his paper, a
thorough revision is needed before his theorem can be said to have
been proven.

There is another difficulty, namely that Brown's theorem does not give
a very direct approach to the Birch and Swinnerton-Dyer conjecture.
This is because it is rather difficult to compute the modular
parameterization and thus the Heegner point, and so the hypotheses of
Brown's theorem are hard to verify.  (The difficulty comes from the
fact that non-archimedean integration seems to be of exponential
complexity in the desired degree of accuracy, in contrast to
archimedean integration which is polynomial time.)  On the other hand
it is quite easy to check whether the $L$-function of $E$ vanishes to
order 0 or 1, these being the only cases where one expects Heegner
points to be of help.  In fact the computation of the entire
$L$-function of $E$ is straightforward and (at least over the rational
function field) can be made efficient using the existence of an
automorphic form corresponding to $E$.  See \cite{RockmoreTan}.  This
situation is the opposite of that in the classical situation; cf.~the
remarks of Birch near the end of $\S4$ of his article in this volume
\cite{Birch}.

In light of this difficulty, a more direct and straightforward
approach to the Birch and Swinnerton-Dyer conjecture for elliptic
curves of rank $\le1$ is called for.  My interest in function field
analogues of Gross-Zagier came about from an effort to understand
Brown's paper and to find a better approach to BSD in this context.

\subsection{Gross-Zagier formulas}
Let us now state what the analogue of the Gross-Zagier formula
\cite{GZ} should be in the
function field context.  Let $E$ be an elliptic curve over $F$ of
conductor $\n\infty$ and with split multiplicative reduction at
$\infty$.  Then for every imaginary quadratic extension $K$ of $F$
satisfying the Heegner hypotheses (namely that every prime dividing
$\n$ is split in $K$), we have a point $P_K\in E(K)$
defined using Heegner points on $X_0(\n)$ and the modular
parameterization.  The desired formula is then
\begin{equation}\label{eq:GZ1}
L'(E/K,1)=a\langle P_K,P_K\rangle
\end{equation}
where $\langle,\rangle$ is the N\'eron-Tate canonical height on $E$ and
$a$ is an explicit non-zero constant.  Because equality of analytic and
algebraic ranks implies the refined BSD conjecture, the exact value of
$a$ is not important for us.

The left hand side of this formula is also a special value of the
$L$-function of an automorphic form (namely, the $f$ such that
$L(E,\chi,s)=L(f,\chi,s)$) and Equation~\ref{eq:GZ1} is a special case
of a more general formula which applies to automorphic forms without
the assumption that their Hecke eigenvalues are integers.  Let $S$ be
the vector space of complex valued cuspidal automorphic forms on
$\GL_2$ over $F$ which have level $\Gamma_0(\n\infty)$, central
character $|\cdot|^{-2}$, and which are special at $\infty$.  (As
discussed in Subsection~\ref{ss:AnalyticMod}, this is the analogue of
$S_2(\Gamma_0(N))$.)  Then we have a Petersson inner product
\begin{equation*}
(\ ,\ ):S\times S\to\C
\end{equation*}
which is positive definite Hermitian.  For $f\in S$, let $L_K(f,s)$ be
the $L$-function of the base change of $f$ to a form on $\GL_2$ over
$K$.  (This form can be shown to exist using a Rankin-Selberg integral
representation and Weil's converse theorem.)  
Then the function $f\mapsto L_K'(f,1)$ is a linear map
$S\to\C$ and so there exists a unique element $h_K\in S$ such that
\begin{equation*}
(f,h_K)=L'_K(f,1)
\end{equation*}
for all $f\in S$.  

For $h\in S$, let $c(h,\m)$ be the $\m$-th Fourier
coefficient of $h$.  Then a formal Hecke algebra argument, as in the
classical case, shows that the desired Gross-Zagier formula
\ref{eq:GZ1} (and its more general version mentioned above) follows
from the following equalities between Fourier coefficients and heights
on $J_0(\n)$:
\begin{equation}\label{eq:GZ2}
c(h_K,\m)=a\langle Q_K,T_\m Q_K\rangle
\end{equation}
for all effective divisors $\m$ prime to $\n\infty$.  Here $T_\m$ is the
Hecke operator on $J_0(\n)$ indexed by $\m$ and $\langle,\rangle$
is the canonical height pairing on $J_0(\n)$.

From now on, by ``Gross-Zagier formula'' we will mean the sequence of
equalities~\ref{eq:GZ2}.

\subsection{R\"uck-Tipp}
R\"uck and Tipp were the first to write down a function field analogue
of the Gross-Zagier formula \cite{RuckTipp}.  They work over $F=\Fq(t)$ with
$q$ odd, and $\infty$ the standard place at infinity $t=\infty$ (so
their $\infty$ has degree 1).  They assume that $\n$ is square free
and that $K=F(\sqrt D)$ where $D$ is an {\it irreducible\/} polynomial
in $\Fq[t]$.  Under these hypotheses, they checked the
equalities~\ref{eq:GZ2} for all $\m$ prime to $\n\infty$, which yields
the formula~\ref{eq:GZ1}.  This gives some instances of the conjecture
of Birch and Swinnerton-Dyer, under very restrictive hypotheses.

Their paper follows the method of Gross and Zagier \cite{GZ} quite
closely (which is not to say that the analogies are always obvious or
easy to implement!). They use the Rankin-Selberg method and a
holomorphic projection operator to compute the Fourier coefficients of
$h_K$.  The height pairing is decomposed as a sum of local terms and,
at finite places, the local pairing is given as an intersection
number, which can be computed by counting isogenies between Drinfeld
modules over a finite field.  The local height pairing at $\infty$ is
also an intersection number and one might hope to use a moduli
interpretation of the points on the fibre at $\infty$ to calculate the
local height.  But to my knowledge, no one knows how to do this.
Instead, R\"uck and Tipp compute the local height pairing using
Green's functions on the Drinfeld upper half plane.  This is a very
analytic way of computing a rational number, but it matches well with
the computations on the analytic side of the formula.

\subsection{P\'al and Longhi}
P\'al and Longhi worked (independently) on function field analogues of
the Bertolini-Darmon \cite{BD1} $p$-adic construction of Heegner
points.  Both work over a general function field $F$ of odd
characteristic.  Let $E$ be an elliptic curve over $F$ with conductor
$\n\infty$ and which is split multiplicative at $\infty$.  Let $K$ be
a quadratic extension in which $\infty$ is inert and which satisfies
the Heegner hypotheses with respect to $E$.  Also let $H_n$ be the
ring class field of $K$ of conductor $\infty^n$ and set
$G=\invlim\gal(H_n/K)$.

P\'al \cite{Pal} used ``Gross-Heegner'' points, as in Bertolini-Darmon
(following Gross~\cite{GrossM}), to construct an element $\LL(E/K)$ in
the completed group ring $\Z[[G]]$ which interpolates suitably
normalized special values $L(E/K,\chi,1)$ for finite order characters
$\chi$ of $G$.  It turns out that $\LL(E/K)$ lies in the augmentation
ideal $I$ of $\Z[[G]]$ and so defines an element $\LL'(E/K)$ in
$I/I^2\cong\O_{K_\infty}^\times/\O_{F_\infty}^\times
\cong\O_{K_\infty,1}^\times$.  (Here $F_\infty$ and $K_\infty$ are the
completions at $\infty$ and $\O_{K_\infty,1}^\times$ denotes the
1-units in $\O_{K_\infty}$.)  Since $E$ is split multiplicative at
$\infty$, we have a Tate parameterization $K_\infty^\times\to
E(K_\infty)$ and P\'al shows that the image of $\LL'(E/K)$ in
$E(K_\infty)$ is a global point.  More precisely, if $E$ is a ``strong
Weil curve,'' then P\'al's point  is $P_K-\overline{P}_K$ where $P_K$ is the Heegner
point discussed above and $\overline{P}_K$ is its ``complex
conjugate.''  It follows that if $\LL'(E/K)$ is non-zero, then the
Heegner point is of infinite order and so $\rk_\Z E(K)$ is at least
one.  One interesting difference between P\'al's work and \cite{BD1}
is that in the latter, there are 2 distinguished places, namely
$\infty$, which is related to the classical modular parameterization,
and $p$, which is related to the Tate parameterization.  In P\'al's
work, the role of both of these primes is played by the prime $\infty$
of $F$.  This means that his result is applicable in more situations
than the naive analogy would predict---$E$ need only have split
multiplicative reduction at one place of $F$.

Longhi \cite{Longhi} also gives an $\infty$-adic construction of a
Heegner point.  Whereas P\'al follows \cite{BD1}, Longhi's point of
view is closer to that of \cite{BD2}.  His $\infty$-adic $L$-element
$\LL(E/K)$ is constructed using $\infty$-adic integrals, following the
approach of Schneider~\cite{Schneider} and a multiplicative version of
Teitelbaum's Poisson formula \cite{Teitelbaum}.  Unfortunately, there
is as yet no connection between his $\infty$-adic $\LL(E/K)$ and
special values of $L$-functions.

Both of these works have the advantage of avoiding intricate height
computations on Drinfeld modular curves, as in \cite{GZ}. (P\'al's
work uses heights of the much simpler variety considered in
\cite{GrossM}.)  On the other hand, they do not yet have any direct
application to the conjecture of Birch and Swinnerton-Dyer, because
presently we have no direct link between the $\infty$-adic
$L$-derivative $\LL'(E/K)$ and the classical $L$-derivative
$L'(E/K,1)$.

\subsection{My work on BSD for rank 1}
My interest in this area has been less in analogues of the
Gross-Zagier formula or Kolyvagin's work over function fields {\it per
  se\/}, and more in their applications to the Birch and
Swinnerton-Dyer conjecture itself.  The problem with a raw
Gross-Zagier formula is that it only gives the BSD conjecture with
parasitic hypotheses.  For example, to have a Drinfeld modular
parameterization, and thus Heegner points, the elliptic curve must
have split multiplicative reduction at some place and the existence of
such a place presumably has nothing to do with the truth of the
conjecture.  Recently, I have proven a non-vanishing result which
when combined with a suitable Gross-Zagier formula leads to a clean,
general statement about Birch and Swinnerton-Dyer: ``If $E$ is an
elliptic curve over a function field $F$ of characteristic $>3$ and
$\ord_{s=1}L(E/F,s)\le1$, then the Birch and Swinnerton-Dyer
conjecture holds for $E$.''  In the remainder of this section I will
describe the non-vanishing result, and then give the statement and
status of the Gross-Zagier formula I have in mind.

Thus, let $E$ be an elliptic curve over $F$ with
$\ord_{s=1}L(E/F,s)\le1$; for purposes of BSD we may as well assume
that $\ord_{s=1}L(E/F,s)=1$ and that $E$ is non-isotrivial.  Because
$j(E)\not\in\Fq$, it has a pole at some place of $F$, i.e., $E$ is
potentially multiplicative there.  Certainly we can find a finite
extension $F'$ of $F$ such that $E$ has a place of split
multiplicative reduction and it will suffice to prove BSD for $E$ over
$F'$.  But, to do this with Heegner points, we must be able to choose
$F'$ so that $\ord_{s=1}L(E/F',s)$, which is {\it a priori\/} $\ge1$,
is equal to 1.  This amounts to a non-vanishing statement for a
(possibly non-abelian) twist of $L(E/F,s)$, namely
$L(E/F',s)/L(E/F,s)$.  Having done this, a similar issue comes up in
the application of a Gross-Zagier formula, namely, we must find a
quadratic extension $K/F'$ satisfying the Heegner hypotheses such that
$\ord_{s=1}L(E/K,s)=\ord_{s=1}L(E/F',s)=1$.  This amounts to a
non-vanishing statement for quadratic twists of $L(E/F',s)$ by
characters satisfying certain local conditions.  This issue also
comes up in the applications of the classical Gross-Zagier formula and
is dealt with by automorphic methods.  Recently, I have proven a very
general non-vanishing theorem for motivic $L$-functions over function
fields using algebro-geometric methods which when applied
to elliptic curves yields the following result:

\begin{thm} \cite{UlmerGNV}
Let $E$ be a non-constant elliptic curve over a function field $F$ of
characteristic $p>3$.  Then there exists a finite separable extension
$F'$ of $F$ and a quadratic extension $K$ of $F'$ such that the
following conditions are satisfied:
\begin{enumerate}
\item $E$ is semistable over $F'$, i.e., its conductor is square-free.
\item $E$ has split multiplicative reduction at some place of $F'$
which we call $\infty$.
\item $K/F'$ satisfies the Heegner hypotheses with respect to $E$ and
$\infty$.  In other words, $K/F'$ is split at every place
$v\neq\infty$ dividing the conductor of $E$ and it is not split at
$\infty$. 
\item $\ord_{s=1}L(E/K,s)$ is odd and at most $\ord_{s=1}L(E/F,s)+1$.
In particular, if $\ord_{s=1}L(E/F,s)=1$, then
$\ord_{s=1}L(E/K,s)=\ord_{s=1}L(E/F',s)=1$. 
\end{enumerate}
\end{thm}

This result, plus a suitable Gross-Zagier formula, yields the desired
theorem.  Indeed, by point (2), $E$ admits a Drinfeld modular
parameterization over $F'$ and by point (3) we will have a Heegner
defined over $K$.  Point (4) (plus GZ!) guarantees that the Heegner
point will be non-torsion and so we have $\rk E(K)\ge1$.  As we have
seen, this implies BSD for $E$ over $K$ and thus also over $F$.  Point
(1) is included as it makes the needed GZ formula a little more
tractable.  Also, although it is not stated in the theorem, it is
possible to specify whether the place $\infty$ of $F'$ is inert or
ramified in $K$ and this too can be used to simplifly the Gross-Zagier
calculation.

Thus, the Gross-Zagier formula we need is in the following context:
the base field $F'$ is arbitrary but the level $\n$ is square-free and
we may assume that $\infty$ is inert (or ramified) in $K$.  It would
perhaps be unwise to write too much about a result which is not
completely written and refereed, so I will just say a few words.  The
proof follows closely the strategy of Gross and Zagier, with a few
simplifications due to Zhang~\cite{Zhang}.  One computes the analytic
side of \ref{eq:GZ2} using the Rankin-Selberg method and a holomorphic
projection and the height side is treated using intersection theory at
the finite places and Green's functions at $\infty$.  Because we work
over an arbitrary function field, our proofs are necessarily adelic.
Also, in the analytic part we emphasize the geometric view of
automorphic forms, namely that they are functions on a moduli space of
rank 2 vector bundles on $\Curve$.  The full details will appear in
\cite{UlmerGZ}.

\section{Ranks over function fields}\label{s:ffRanks}
We now move beyond rank 1 and consider the rank conjecture for
elliptic curves over function fields.  Recall from Section~\ref{s:BSD}
the notions of constant, isotrivial, and non-isotrivial for elliptic
curves over function fields.  Our purpose in this section is to explain
constructions of isotrivial and non-isotrivial elliptic curves over
$\Fp(t)$ whose Mordell-Weil groups have arbitrarily large rank.  These
curves turn out to have asymptotically maximal rank, in a sense which
we will explain in Section~\ref{s:bounds}.

\subsection{The Shafarevitch-Tate construction}
First, note that if $E$ is a constant elliptic curve over
$F=\Fq(\Curve)$ based on $E_0$, then $E(F)\cong\mor_{\Fq}(\Curve,E_0)$
(morphisms defined over $\Fq$) and the torsion subgroup of $E(F)$
corresponds to constant morphisms.  Since a morphism $\Curve\to E$ is
determined up to translation by the induced map of Jacobians, we have
$E(F)/tor\cong\Hom_{\Fq}(J(\Curve),E)$ where $J(\Curve)$ denotes the
Jacobian of $\Curve$.

The idea of Shafarevitch and Tate \cite{TateShaf} was to take $E_0$ to
be supersingular and to find a curve $\Curve$ over $\Fp$ which is
hyperelliptic and such that $J(\Curve)$ has a large number of 
factors isogenous to $E_0$.  If $E$ denotes the constant curve over
$\Fp(t)$ based on $E_0$, then it is clear that $E(\Fp(t))$ has rank 0.
On the other hand, over the quadratic extension $F=\Fp(\Curve)$,
$E(F)/tor\cong\Hom_{\Fq}(J(\Curve),E_0)$ has large rank.  Thus if we let
$E'$ be the twist of $E$ by the quadratic extension $F/\Fq(t)$, then
$E'(\Fq(t))$ has large rank.  Note that $E'$ is visibly isotrivial.

To find such curves $\Curve$, Tate and Shafarevitch considered
quotients of the Fermat curve of degree $p^n+1$ with $n$ odd.  The
zeta functions of Fermat curves can be computed in terms of Gauss
sums, and in the case of degree of the form $p^n+1$, the relevant
Gauss sums are easy to make explicit.  This allows one to show that
the Jacobian is isogenous to a product of supersingular elliptic
curves over $\Fpbar$ and has a supersingular elliptic curve as isogeny
factor to high multiplicity over $\Fp$.

We remark that the number of factors of $J(\Curve)$ which are
isogenous to $E_0$ may go up under extension of the ground field, and
so the rank of $E'$ may also go up.  In fact, the rank of the
Shafarevitch-Tate curves goes up considerably: if the rank of
$E'(\Fp(t))$ is $r$, then the rank of $E'(\Fpbar(t))$ is of the order
$2\log_p(r)r$.

It has been suggested by Rubin and Silverberg that one might be able
to carry out a similar construction over $\Q(t)$, i.e., one might try
to find hyperelliptic curves $\Curve$ defined over $\Q$ whose
Jacobians have as isogeny factor a large number of copies of some
elliptic curve.  The obvious analogue of the construction above would
then produce elliptic curves over $\Q(t)$ of large rank.  In
\cite{RubinSilverberg} they use this idea to find many elliptic curves
of rank $\ge3$.  Unfortunately, it is not at all evident that curves
$\Curve$ such that $J(\Curve)$ has an elliptic isogeny factor to high
multiplicity exist, even over $\C$.

Back to the function field case: We note that isotrivial elliptic
curves are very special and seem to have no analogue over $\Q$.  Thus
the relevance of the Shafarevitch-Tate construction to the rank
question over $\Q$ is not clear.  In the next subsection we explain
a construction of {\it non-isotrivial\/} elliptic curves over $\Fp(t)$
of arbitrarily large rank.

\subsection{Non-isotrivial elliptic curves of large rank}
\label{ss:HighRankCurves}
In \cite{Shioda}, Shioda showed that one could often compute the
Picard number of a surface which is dominated by a Fermat surface.  He
applied this to write down elliptic curves over $\Fpbar(t)$ (with
$p\equiv3\mod4$) of arbitrarily large rank, using supersingular Fermat
surfaces (i.e., those whose degrees divide $p^n+1$ for some $n$).  I
was able to use the idea of looking at quotients of Fermat surfaces
and a different method of computing the rank to show the existence of
elliptic curves over $\Fp(t)$ (any $p$) with arbitrarily large rank.
Here is the precise statement:

\begin{thm}\label{thm:ranks}\cite{UlmerR}
Let $p$ be a prime, $n$ a positive integer, and $d$ a divisor of
$p^n+1$.  Let $q$ be a power of $p$ and let $E$ be the elliptic curve
over $\Fq(t)$ defined by
\begin{equation*}
y^2+xy=x^3-t^d.
\end{equation*}
Then the $j$-invariant of $E$ is not in $\Fq$, the conjecture of Birch
and Swinnerton-Dyer holds for $E$, and the rank of $E(\Fq(t))$ is
\begin{equation*}
\sum_{\substack{e|d\\e\nodiv6}}\frac{\phi(e)}{o_e(q)}
+\begin{cases}
        0&\text{if $2\nodiv d$ or $4\nodiv q-1$}\\
        1&\text{if $2|d$ and $4|q-1$}
\end{cases}
+\begin{cases}
        0&\text{if $3\nodiv d$}\\
        1&\text{if $3|d$ and $3\nodiv q-1$}\\
        2&\text{if $3|d$ and $3|q-1$.}
\end{cases}
\end{equation*}
Here $\phi(e)$ is the cardinality of $(\Z/e\Z)^\times$ and $o_e(q)$
is the order of $q$ in $(\Z/e\Z)^\times$.
\end{thm}

In particular, if we take $d=p^n+1$ and $q=p$, then the rank of $E$
over $\Fp(t)$ is at least $(p^n-1)/2n$.  On the other hand, if we take
$d=p^n+1$ and $q$ to be a power of $p^{2n}$, then the rank of $E$ over
$\Fq(t)$ is $d-1=p^n$ if $6\nodiv d$ and $d-3=p^n-2$ if $6|d$.  Note
that the rank may increase significantly after extension of $\Fq$.

Here is a sketch of the proof: by old work of Artin and Tate
\cite{TateBSD}, the conjecture of Birch and Swinnerton-Dyer of $E$ is
equivalent to the Tate conjecture for the elliptic surface
$\EE\to\P^1$ over $\Fq$ attached to $E$.  (The relevant Tate
conjecture is that $-\ord_{s=1}\zeta(\EE,s)=\rk_\Z NS(\EE)$ where
$NS(\EE)$ denotes the N\'eron-Severi group of $\EE$.)  The equation of
$E$ was chosen so that $\EE$ is dominated by the Fermat surface {\it
  of the same degree\/} $d$.  (The fact that the equation of $E$ has 4
monomials is essentially enough to guarantee that $\EE$ is dominated
by some Fermat surface; getting the degree right requires more.)
Since the Tate conjecture is known for Fermat surfaces, this implies
it also for $\EE$ (and thus BSD for $E$).  Next, a detailed analysis
of the geometry of the rational map $F_d\ratto\EE$ allows one to
calculate the zeta function of $\EE$ in terms of that of $F_d$, i.e.,
in terms of Gauss and Jacobi sums.  Finally, because $d$ is a divisor
of $p^n+1$, the relevant Gauss sums are all supersingular (as in the
Shafarevitch-Tate case) and can be made explicit.  This gives the
order of pole of $\zeta(\EE,s)$ at $s=1$ and thus the order of zero of
$L(E/\Fq(t),s)$ at $s=1$, and thus the rank.

We note that the proof does not explicitly construct any points,
although it does suggest a method to do so.  Namely, using the large
automorphism group of the Fermat surface, one can write down curves
which span $NS(F_d)$ and use these and the geometry of the map
$F_d\ratto\EE$ to get a spanning set for $NS(\EE)$ and thus a spanning
set for $E(\Fq(t))$.  It looks like an interesting problem to make
this explicit, and to consider the heights of generators of
$E(\Fq(t))$ and its Mordell-Weil lattice.

\subsection{Another approach to high rank curves}
The two main parts of the argument of
Subsection~\ref{ss:HighRankCurves} could be summarized as follows: (i)
one can deduce the Tate conjecture for $\EE$ and thus the BSD
conjecture for $E$ from the existence of a dominant rational map from
the Fermat surface $F_d$ to the elliptic surface $\EE$ attached to
$E$; and (ii) a detailed analysis of the geometry of the map
$F_d\ratto\EE$ allows one to compute the zeta function of $\EE$ and
thus the $L$-function of $E$, showing that it has a large order zero
at $s=1$.  

Ideas of Darmon give an alternative approach to the second part of
this argument (showing that the $L$-function has a large order zero at
$s=1$) and may lead (subject to further development of Gross-Zagier
formulas in the function field case) to an alternative approach to the
first part of the argument (the proof of BSD).  Darmon's idea is quite
general and leads to the construction of many elliptic curves over
function fields of large rank (more precisely, provably of large
analytic rank and conjecturally of large algebraic rank.)  Here we
will treat only the special case of the curve considered in
Subsection~\ref{ss:HighRankCurves} and we refer to his article in this
volume \cite{Darmon} for details of the general picture.

Let $q=p^n$ ($p$ any prime), $d=q+1$, and define $F=\Fq(u)$,
$K=\Fqq(u)$, and $H=\Fqq(t)$ where $u=t^d$.  Then $H$ is Galois over
$F$ with dihedral Galois group.  Indeed $\gal(H/K)$ is cyclic of order
$d$ and because $q\equiv-1\mod d$, the non-trivial element of
$\gal(K/F)\cong\gal(\Fqq/\Fq)$ acts on $\gal(H/K)$ by inversion.  Let
$E$ be the elliptic curve over $F$ defined by the equation
\begin{equation*}
y^2+xy=x^3-u.
\end{equation*}
Over $H$, this is the curve discussed in Subsection~\ref{ss:HighRankCurves}.  

The $L$-function of $E$ over $H$ factors into a product of twisted
$L$-functions over $K$:
\begin{equation*}
L(E/H,s)=\prod_{\chi\in\hat G}L(E/K,\chi,s)
\end{equation*}
where the product is over the $d$ characters of $G=\gal(H/K)$.
Because $H/F$ is a dihedral extension and $E$ is defined over $F$, we
have the equality $L(E/K,\chi,s)=L(E/K,\chi^{-1},s)$. Thus the
functional equation
\begin{align*}
L(E/K,\chi,s)&=W(E/K,\chi,)q^{sd_{E,\chi}}L(E/K,\chi^{-1},2-s)\cr
&=W(E/K,\chi)q^{sd_{E,\chi}}L(E/K,\chi,2-s)
\end{align*}
(where $W(E/K,\chi)$ is the ``root number'' and $d_{E,\chi}$ is the
degree of $L(E/K,\chi,s)$ as a polynomial in $q^{-2s}$) may force a
zero of $L(E/K,\chi,s)$ at the critical point $s=1$.  This is indeed
what happens: A careful analysis shows that $W(E/K,\chi)$ is $+1$ if
$\chi$ is trivial or of order exactly 6 and it is $-1$ in all other
cases.  Along the way, one also finds that $d_{E,\chi}$ is 0 if $\chi$
is trivial or of order exactly 6 and is 1 in all other cases.  Thus
$L(E/K,\chi,s)$ is equal to 1 if $\chi$ is trivial or of order exactly
6 and is equal to $(1-q^{-2s})$ and vanishes to order 1 at $s=1$ if
$\chi$ is non-trivial and not of order exactly 6.  We conclude that
$\ord_{s=1}L(E/H,s)$ is $d-3$ if 6 divides $d$ and $d-1$ if not.

Of course one also wants to compute the $L$-function of $E$ over
$H_0=\Fp(t)$.  In this case, the $L$-function again factors into a
product of twists, but the twists are by certain, generally
non-abelian, representations of the Galois group of the Galois closure
of $\Fp(t)$ over $\Fp(u)$.  (The Galois closure is $H$ and the Galois
group is the semidirect product of $\gal(H/K)$ with
$\gal(\Fqq(u)/\Fp(u))$.  See \cite[\S3]{UlmerGNV} for more on this
type of situation.)  We will not go into the details, but simply note
that in order to compute the $L$-function $L(E/H_0,s)$ along the lines
above, one needs to know the root numbers $W(E/\Fr(u),\chi)$ where
$r=p^{o_e}$, $o_e$ is the order of $p$ in $\Z/e\Z$, and $e$ is the
order of the character $\chi$. It turns out that each of the twisted
$L$-functions has a simple zero at $s=1$.

This calculation of $L(E/H,s)$ and $L(E/\Fp(t),s)$ seems to be of
roughly the same difficulty as the geometric one in \cite{UlmerR}
because the ``careful analysis" of the root numbers $W(E/K,\chi)$ and
$W(E/\Fr(u),\chi)$ is somewhat involved, especially if one wants to
include the cases $p=2$ or 3. (I have only checked that the answer
agrees with that in \cite{UlmerR} when $p>3$.)  Calculation of the
root numbers requires knowing the local representations of
decompostion groups on the Tate module at places of bad reduction and
eventually boils down to analysing some Gauss sums.  The
Shafarevitch-Tate lemma on supersingular Gauss sums (Lemma~8.3 of
\cite{UlmerR}) is a key ingredient.

Regarding the problem of verifying the BSD conjecture for $E/H$, note
that $K/F$ may be viewed as an ``imaginary" quadratic extension, and
$H/K$ is the ring class extension of conductor $\n=(0)(\infty)$.
Because most of the twisted $L$-functions $L(E/K,\chi,s)$ vanish
simply, we might expect to constuct points in $(E(H)\tensor\C)^\chi$
using Heegner points and show that they are non-trivial using a
Gross-Zagier formula.  But the relevant Gross-Zagier formula here
would involve Shimura curve analogs of Drinfeld modular curves (since
the extension $K/F$ does not satisfy the usual Heegner hypotheses) and
such a formula remains to be proven.  Perhaps Darmon's construction
will provide some motivation for the brave soul who decides to take on
the Gross-Zagier formula in this context!  On the other hand, Darmon's
paper has examples of curves where Heegner points on standard Drinfeld
modular curves should be enough to produce high rank elliptic curves
over $\Fp(t)$.

\numberwithin{equation}{section}
\section{Rank bounds}\label{s:bounds}
We now return to a general function field $F=\Fq(\Curve)$ and a
general non-isotrivial elliptic curve $E$ over $F$.  Recall that the
conductor $\n$ of $E$ is an effective divisor on $\Curve$ which is
supported precisely at the places where $E$ has bad reduction.  

It is natural to ask how quickly the ranks of elliptic curves over $F$
can grow in terms of their conductors.  As discussed in
Section~\ref{s:BSD}, we have the inequality
\begin{equation*}
\rk_\Z E(F)\le\ord_{s=1}L(E/F,s).
\end{equation*}
Also, one knows
that that $L(E/F,s)$ is a polynomial in $q^{-s}$ of degree
$4g-4+\deg\n$ where $g$ is the genus of $\Curve$.  (This comes from
Grothendieck's cohomological expression for the $L$-function and the
Grothendieck-Ogg-Shafarevitch Euler characteristic formula.)  Thus we
have a bound
\begin{equation}
\rk_\Z E(F)\le \ord_{s=1}L(E/F,s)\le 4g-4+\deg\n.
\end{equation}

This bound is geometric in the sense that it does not involve the size
of the finite field $\Fq$; the same bound holds for $\rk_\Z
E(\Fqbar(\Curve))$.  On the other hand, as we have seen above, the
rank can change significantly after extension of $\Fq$.  It is thus
natural to ask for a more arithmetic bound, i.e., one which is
sensitive to $q$.

Such a bound was proven by Brumer~\cite{Brumer}, using Weil's
``explicit formula'' technique, along the lines of Mestre's bound for
the rank of an elliptic curve over $\Q$.  Brumer's result is
\begin{equation}
\rk_\Z E(F)\le \ord_{s=1}L(E/F,s)\le\frac{4g-4+\deg(\n)}{2\log_q\deg(\n)}
+C\frac{\deg(\n)}{(\log_q\deg(\n))^2}
\end{equation}

Note that this bound is visibly sensitive to $q$ and is an improvement
on the geometric bound when $\deg\n$ is large compared to $q$.

Here is a sketch of Brumer's proof: let $\Lambda(s)=q^{Ds/2}L(E/F,s)$
where $D=4g-4+\deg\n$.  Then $\Lambda(s)$ is a Laurent polynomial in
$q^{-s/2}$ and so is periodic in $s$ with period $4\pi i/\ln q$;
moreover, we have the functional equation
$\Lambda(s)=\pm\Lambda(2-s)$.  Our task is to estimate the order of
vanishing at $s=1$ of $\Lambda$ or equivalently, the residue at $s=1$
of the logarithmic derivative $\Lambda'/\Lambda$ with respect to $s$.
Let us consider the line integral
\begin{equation*}
I=\oint\Phi d\log\Lambda=\oint\Phi\frac{\Lambda'}{\Lambda}\,ds
\end{equation*}
where $\Phi$ is a suitable test function to be chosen later and the
contour of integration is
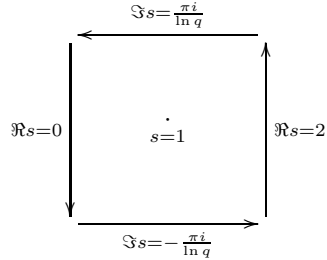
\begin{equation*}
\xymatrix{\ar[dd]_{\Re s=0}&&\ar[ll]_{\Im s=\frac{\pi i}{\ln q}}\\
&\underset{s=1}{\cdot}\\
\ar[rr]_{\Im s=-\frac{\pi i}{\ln q}}&&\ar[uu]_{\Re s=2}}
\end{equation*}
(Note that $d\log\Lambda$ is periodic with period $2\pi i/\ln q$.
Also, we would have to shift the contour slightly if $L(E/F,s)$ has a
zero at $1\pm\pi i/\ln q$.)  We assume that $\Phi(s)$ is non-negative
on the line $\Re s=1$.  By the Riemann hypothesis for $L(E/F,s)$, all
the zeroes of $\Lambda(s)$ lie on this line and so
\begin{equation*}
\Phi(1)\ord_{s=1}L(E/F,s)=\Phi(1)\res_{s=1}\frac{\Lambda'}\Lambda
\le \sum_{\rho}\Phi(\rho)=I
\end{equation*}
where $\rho$ runs over the zeroes of $L(E/F,s)$ inside the contour of
integration counted with multiplicities.  Now we assume in addition that
$\Phi$ is a Laurent polynomial in $q^{-s}$ (so periodic with period
$2\pi i/\ln q$) and that it satisfies the functional equation
$\Phi(s)=\Phi(2-s)$.  Using the functional equation and periodicity of
the integrand, the integral $I$ is equal to
\begin{equation}
2\int_{2-\pi i/\ln q}^{2+\pi i/\ln q}\Phi\frac{\Lambda'}{\Lambda}\,ds.
\end{equation}
Now the integration takes place entirely in the region of convergence
of the Euler product defining $L(E/F,s)$ and so we can expand the
integrand in a series and estimate the terms using the Riemann
hypothesis for curves over finite fields.  Finally, Brumer makes a
clever choice of test function $\Phi$ which yields the desired
estimate.  (More precisely, he considers a sequence of test functions
satisfying the hypotheses which when restricted to $\Re s=1$ converge
to the Dirac delta function at $s=1$---up to a change of variable,
this is essentially the Fej\'er kernel---and then chooses $\Phi$ to be
a suitable element of this sequence.)

Note the strongly analytic character of this proof.  For example, it
does not use the fact that there is massive cancellation in the series
for $L(E/F,s)$ so that the $L$-function is really a polynomial in
$q^{-s}$ of degree $D$!

Let $E_n$ be the curve of Theorem~\ref{thm:ranks} with $d=p^n+1$.
Then it turns out that the degree of the conductor of $E_n$ is $p^n+2$
if $6|d$ and $p^n+4$ if $6\nodiv d$.  One sees immediately that the
geometric bound is sharp when $q$ is a power of $p^{2n}$ and the main
term of the arithmetic bound is met when $q=p$.  Thus both the
geometric and arithmetic bounds give excellent control on ranks.

The rest of this paper is devoted to considering various questions
which arise naturally by analogy from the existence and sharpness of
these two types of bounds, geometric and arithmetic.  

\section{Ranks over number fields}\label{s:nfRanks}
We now turn to analogous situations, starting with the case where the
ground field $F$ is either $\Q$ or a number field.  Throughout, we
assume that the $L$-series of elliptic curves have good analytic
properties, namely analytic continuation, boundedness in vertical
strips, and the standard functional equation.  (This is of course now
known for elliptic curves over $\Q$ by the work of Wiles and others,
but is still open for a general number field $F$.) We also assume the
conjecture of Birch and Swinnerton-Dyer so that ``rank'' can be taken
to mean either analytic rank ($\ord_{s=1}L(E/F,s)$) or algebraic rank
($\rk_\Z E(F)$); alternatively the reader may interpret each question
or conjecture involving an unqualified ``rank'' to be two statements,
one about analytic rank, the other about algebraic rank.

The Brumer bound discussed in the last section was modeled on work of
Mestre \cite{Mestre}, who proved, along lines quite similar to those
sketched above, a bound on analytic ranks of the following form:
\begin{equation}\label{eqn:Mestre}
\ord_{s=1}L(E/\Q,s)=O\left(\frac{\log N}{\log\log N}\right)
\end{equation}
where $E$ is an elliptic curve over $\Q$ of conductor $N$.  To see the
analogy, note that the degree functon on divisors is a kind of
logarithm and so $\deg\n$ is an analogue of $\log N$.  To obtain this
bound, Mestre assumes the Generalized Riemann Hypothesis for
$L(E/\Q,s)$ and he actually proves a more general statement about
orders of vanishing for $L$-series of modular forms.  Assuming good
analytic properties and the generalized Riemann hypothesis, his
argument extends readily to elliptic curves over number fields; in
this case $N$ should be replaced with the norm from $F$ to $\Q$ of the
conductor of $E$ times the absolute value of the discriminant of $F$. 

There is some evidence that the Mestre bound should be asymptotically
sharp.  First of all, it gives excellent results for small $N$.
Secondly, in the function field case, the analogous bound is sharp.
Moreover, the proof of the bound in the function field case does not
use strongly any special features of that situation, such as the fact
that the $L$-function is really a polynomial.  Motivated by these
facts, I make the following conjecture about the sharpness of the
Mestre bound.

\begin{conj}  Fix a number field $F$ and for each positive integer
  $N$, define $r_F(N)$ by
\begin{equation*}
r_F(N)=\max\{\rk_\Z(E(F)) | E/F\text{ with
}\text{Norm}_{F/\Q}(\n_E)=N\}
\end{equation*}
where the maximum is taken over all elliptic curves $E$ over $F$ with
conductor $\n_E$ satisfying $\text{Norm}_{F/\Q}(\n)=N$;
if there are no such curves, we set $r_F(N)=0$.  Then we have
\begin{equation*}
\limsup_N\frac{r_F(N)}{\log N/\log\log N}>0
\end{equation*}
\end{conj}

By the generalization of the Mestre bound
to number fields, the limit in the conjecture is finite.

If $E$ is an elliptic curve over $\Q$, let $N_\Q(E)$ be its conductor
and let $N_F(E)$ be the norm from $F$ to $\Q$ of the conductor of $E$
viewed as elliptic curve over $F$.  Then there is a constant
$C$ depending only on $F$ such that
\begin{equation*}
1\le\frac{N_\Q(E)^{[F:\Q]}}{N_F(E)}\le C
\end{equation*}
for all elliptic curves $E$ over $\Q$.  Indeed, if $N$ is an integer,
then $\text{Norm}_{F/\Q}(N)=N^{[F:\Q]}$ and since the conductor of $E$
  over $F$ is a divisor of $N_\Q(E)$ (viewed as an ideal of $F$), we
  have $1\le{N_\Q(E)^{[F:\Q]}}/{N_F(E)}$.  Since $N_\Q(F)$ divided by
  the conductor of $E$ over $F$ is divisible only by ramified primes
  and these occur with bounded exponents~\cite{BrumerKramer}, there is
  a constant $C$ such that ${N_\Q(E)^{[F:\Q]}}/{N_F(E)}\le C$.  These
  inequalities show that a sequence of elliptic curves proving the
  conjecture over $\Q$ also proves the conjecture for a general number
  field $F$.

Finally, let us remark that there are experts who are skeptical about
this conjecture.  Certain probabilistic models predict that the
denominator should be replaced by its square root, i.e., that the
correct bound is
\begin{equation*}
\ord_{s=1}L(E/\Q,s)\mathrel{\mathop=^?}
O\left(\left(\frac{\log N}{\log\log N}\right)^{1/2}\right).
\end{equation*}
On the other hand, certain random matrix models suggest that the
Mestre bound is sharp.  See the list of problems for the workshop on
random matrices and $L$-functions at AIM, May 2001 ({\tt
http://aimath.org}) for more on this question.

\section{Algebraic rank bounds}\label{s:algRankBounds}
The Mestre and Brumer bounds are analytic in both statement and
proof.  It is interesting to ask whether they can be made more
algebraic.   For example, the Brumer bound is equivalent to a
statement about the possible multiplicity of $q$ as an eigenvalue of
Frobenius on $H^1(\Curve,\FF)$ for a suitable sheaf $\FF$, namely
$R^1\pi_*\Ql$ where $\pi:\EE\to\Curve$ is the elliptic surface
attached to $E/F$.  It seems that statements like this might admit
more algebraic proofs.

There is one situation where such algebraic proofs are available.
Namely, consider an elliptic curve $E$ over a number field or a 
function field $F$ such that $E$ has an $F$-rational 2-isogeny $\phi:E\to
E'$.  Then the Selmer group for multiplication by $2$ sits in an exact
sequence
\begin{equation*}
\text{Sel}(\phi)\to\text{Sel}(2)\to\text{Sel}(\check\phi)
\end{equation*}
where $\check\phi:E'\to E$ is the dual isogeny.  The orders of the
groups $\text{Sel}(\phi)$ and $\text{Sel}(\check\phi)$ can be crudely
and easily estimated in terms of $\omega(N)$, the number of primes dividing the
conductor $N$ of $E$, and a constant depending only on $F$ which involves
the size of its class group and unit group.  This yields a bound
on the rank of the form
\begin{equation*}
\rk_\Z E(F)\le C+2\omega(N).
\end{equation*}
Note that this bound deserves to be called arithmetic because, for
example, in the function field case $F=\Fq(\Curve)$, $\omega(N)$ is
sensitive to $\Fq$ since primes dividing $N$ may split after extension of
$\Fq$.  Note also that it is compatible with the Mestre and Brumer
bounds, since $\omega(N)=O(\log N/\log\log N)$
\cite[p.~355]{HardyWright} in the number field case and
$\omega(N)=O(\deg N/\log\deg N)$ in the function field case.

It is tempting to guess that a similar bound (i.e., $\rk
E(F)=O(\omega(N))$) might be true in general, but there are several
reason for skepticism.  First of all, the estimation of the Selmer
group above breaks down when there is no $F$-rational 2-isogeny.  In
this case, one usually passes to an extension field $F'$ where such an
isogeny exists.  But then the ``constant'' $C$ involves the units and
class groups of $F'$ and these vary with $E$ since $F'$ does.  Given
our current state of knowledge about the size of class groups, the
bounds we obtain are not as good as the Mestre/Brumer
bounds.  This suggests that what is needed is a way to calculate or at
least estimate the size of a Selmer group $\text{Sel}(\ell)$ without
passing to an extension where the multiplication by $\ell$ isogeny
factors.

The second reason for skepticism is that such a bound would imply, for
example, that there is a universal bound on the ranks of elliptic
curves over $\Q$ of prime conductor.  Although we have little
information on the set of such curves (for example, it is not even
known that this set is infinite), the experts seem to be skeptical
about the existence of such a bound.  One fact is that there is an
elliptic curve over $\Q$ with prime conductor and rank 10
\cite{Mestre}, and so the constant in a bound of $O(\omega(N))$ would
have to be at least 10, which does not seem very plausible.  Also, in
\cite{BrumerSilverman}, Brumer and Silverman make a conjecture which
contradicts an $O(\omega(N))$ bound---their conjecture implies that
there should be elliptic curves with conductor divisible only by 2, 3,
and one other prime and with arbitrarily large rank.  There is no
substantial evidence one way or the other for their conjecture, so
some caution is necessary.

Lastly, wild ramification may have some role to play.  Indeed, for
$p=2 $ or $3$ the curves of Section~\ref{s:ffRanks} have conductor
which is divisible only by two primes ($t=0$ and $t=\infty$) and yet
their ranks are unbounded.

Despite all these reasons for skepticism about a bound of the form
$\rk_\Z E(F)\le O(\omega(N))$, it is interesting to ask about the
possibility of estimating ranks or Selmer groups directly, i.e.,
without reducing to isogenies of prime degree.  It seems to me that
there is some hope of doing this in the function field case, at least
in the simplest context of a semistable elliptic curve over the
rational function field.  In this case, ideas from \'etale cohomology
(e.g., \cite[pp.~211-214]{MilneEC}) allow one to express a cohomology
group closely related to the Selmer group as a product of local factors
where the factors are indexed by the places of bad reduction of the
elliptic curve.  

Another approach over function fields is via $p$-descent.  In this
case, there is always a rational $p$-isogeny, namely Frobenius, but,
in contrast to an $\ell$-descent, the places of (good) supersingular
reduction play a role more like places of bad reduction for
$\ell$-descents.  This means that the output of a $p$-descent does not
{\it a priori\/} give a bound for the rank purely in terms of
the conductor and invariants of the ground field.  More work will be
required here to yield interesting results.  See~\cite{Voloch} and
\cite{Ulmerpd} for foundational work on $p$-descents in characteristic
$p$.

\section{Arithmetic and geometric bounds I: cyclotomic fields}
\label{s:cycRanks} \def\Qpc{\Q^{p\text{-cyc}}} 
We now turn to some questions motivated by the existence of a {\it
  pair\/} of bounds, one geometric, one arithmetic.  Let
$\Q_n\subset\Q(\mu_{p^{n+1}})$ be the subfield with $\gal(\Q_n/\Q)$
equal to $\Z/p^n\Z$ and set $\Qpc=\cup_{n\ge0}\Q_n$.  This is the
cyclotomic $\Z_p$-extension of $\Q$.  As is well-known, the extension
$\Qpc/\Q$ may be thought of as an analogue of the extension
$\Fqbar(\Curve)/\Fq(\Curve)$, and this analogy, noted by
Weil~\cite[p.~298]{WeilW}, was developed by Iwasawa into a very
fruitful branch of modern number theory.  There has also been some
traffic in the other direction, e.g., \cite{MazurWiles}.  Let us
consider the rank bounds of Section~\ref{s:bounds} in this light.

Mazur~\cite{Mazur}, in analogy with Iwasawa's work, asked about the
behavior of the Mordell-Weil and Tate-Shafarevitch groups of an
elliptic curve (or abelian variety) defined over $\Q$ as one ascends
the cyclotomic tower.  For example, he conjectured that if $E$ is an
elliptic curve with good, ordinary reduction at $p$, then $E(\Qpc)$
should be finitely generated.  This turns out to be equivalent to the
assertion that $\rk_\Z E(\Q_n)$ is bounded as $n\to\infty$, i.e., it
stabilizes at some finite $n$.

Today, by work of Rohrlich~\cite{Rohrlich}, Kato~\cite{Kato},
Rubin~\cite{Rubin}, and others, this is known to hold even without
the assumption that $E$ has ordinary reduction at $p$.  (But we do
continue to assume that $E$ has good reduction, i.e., that $p$ does
not divide the conductor of $E$.)

Rohrlich proved the analytic version of this assertion, namely that
the analytic rank $\ord_{s=1}L(E/\Q_n,s)$ is bounded as $n\to\infty$.
(Rohrlich's paper is actually about the $L$-functions of modular
forms, but by the work of Wiles and his school, it applies to elliptic
curves.)  Note that
\begin{equation*}
L(E/\Q_n,s)=\prod_\chi L(E/\Q,\chi,s)
\end{equation*}
where $\chi$ ranges over characters of $\gal(\Q_n/\Q)$.  So Rohrlich's
theorem is that for any finite order character $\chi$ of
$\gal(\Qpc/\Q)$ of sufficiently high conductor, $L(E/\Q,\chi,1)\neq0$.
He proves this by considering the average of special values for
conjugate characters $L(E/\Q,\chi^\sigma,1)$ as $\sigma$ varies over a
suitable Galois group and showing that this
average tends to 1 as the conductor of $\chi$ goes to infinity.  Since
$L(E/\Q,\chi^\sigma,1)\neq0$ if and only if $L(E/\Q,\chi,1)\neq0$, this
implies $L(E/\Q,\chi,1)\neq0$.

Work of Rubin, Rubin-Wiles, and Coates-Wiles in the CM case and work
of Kato in the non-CM case (see~\cite[\S8.1]{Rubin} and the references
there) allows us to translate this analytic result into an algebraic
result.  Namely, these authors show that $L(E/\Q,\chi,1)\neq0$ implies
that $(E(\Q_n)\tensor\C)^\chi=0$ where $\chi$ is a character of
$\gal(\Q_n/\Q)$.  This, together with Rohrlich's theorem implies that
the rank of $E(\Q_n)$ stabilizes for large $n$.

Thus for an elliptic curve $E$ over $\Q$ with good reduction at $p$,
$\rk_\Z E(\Qpc)$ is finite and we may ask for a bound.  Since there
are only finitely many $E$ of a given conductor, there is a bound
purely in terms of $p$ and $N$.  The question then is what is the
shape of this bound.  For a fixed $p$, the results of
Section~\ref{s:bounds} might lead one to guess that $\rk_\Z
E(\Qpc)=O(\log N)$ (where the constant of course depends on $p$), but
this is nothing more than a guess.

Rohrlich mentions briefly the issue of an effective bound for the
smallest $q$ such that $L(E/\Q,\chi,1)\neq0$ for all $\chi$ of
conductor $p^n>q$.  He obtains a bound of the form $q=CN^{170}$.
Combined with the Mestre bound~\ref{eqn:Mestre}, this implies that
$\ord_{s=1}L(E/\Q_n,s)$ is bounded for all $n$ by a polynomial in $N$
(which of course depends on $p$).  This bound has recently been
improved by Chinta~\cite{Chinta}.  He points out that
his Theorem~3 (or his Proposition~1 combined with Rohrlich's
arguments) implies the following: If $p$ is an odd prime where $E$ has
good reduction, then for every $\epsilon>0$ there exist
constants $C_\epsilon$ and $e_\epsilon$ such that
\begin{equation*}
\ord_{s=1}L(E/\Q_n,s)\le C_\epsilon p^{e_\epsilon} N^{1+\epsilon}
\end{equation*}
for all $n$.  The exponent $e_\epsilon$ may be taken to be linear in
$1/\epsilon$.  This is of course a weaker bound than the guess $O(\log
N)$; it might be interesting to try to establish the stronger bound on
average.  We remark that Chinta also shows the remarkable result that
there exists an $n_0$ depending on $E$ but {\it independent of $p$\/}
such that $L(E/\Q,\chi,1)\neq0$ for all $\chi$ of conductor $p^n$,
$n>n_0$.

\section{Arithmetic and geometric bounds II: function fields over
  number fields}
\label{s:higherdffRanks} 
Let $K$ be a number field and $\Curve$ a smooth, proper, geometrically
connected curve over $K$.  Let $E$ be a non-isotrivial elliptic curve
over $F=K(\Curve)$ (i.e., $j(E)\not\in K$).  It is known that $E(F)$
is finitely generated \cite{Neron}.

This finite generation, as well as a bound on the rank, can be
obtained by considering the elliptic surface $\pi:\EE\to\Curve$
attached to $E/F$.  As in Section~\ref{s:BSD}, $\EE$ is the unique
elliptic surface over $\Curve$ which is smooth and proper over $K$,
with $\pi$ flat, relatively minimal, and with generic fiber $E/F$.  There is
a close conection between the Mordell-Weil group $E(F)$ and the
N\'eron-Severi group $NS(\EE)$.  Using this, the cycle class map
$NS(\EE)\to H^2(\EE\times\overline{K},\Ql)$ and an Euler
characteristic formula, one obtains the same bound as in the positive
characteristic case, namely:
\begin{equation}\label{eq:geo-bound}
\rk_\Z E(F)\le 4g-4+\deg\n
\end{equation}
where $g$ is the genus of $\Curve$ and $\n$ is the conductor of $E$.

This bound is geometric in that the number field $K$ does not appear
on the right hand side.  
In particular, the bound continues to hold if
we replace $K$ by $\overline K$:
\begin{equation*}
\rk_\Z E(\overline{K}(\Curve))\le 4g-4+\deg\n.
\end{equation*}
Using Hodge theory, this bound can be improved to $4g-4+\deg\n-2p_g$
where $p_g$ is the geometric genus of $\EE$, but this is again a
geometric bound.  It is reasonable to ask if there is a more arithmetic
bound, improving \ref{eq:geo-bound}.

There is some evidence that such a bound exists.
Silverman~\cite{SilvermanES} considered the following situation: Let
$E$ be an elliptic curve over $F=K(t)$ and define $N^*(E)$ to be the
degree of the part of the conductor of $E$ which is prime to $0$ and
$\infty$.  Alternatively, $N^*(E)$ is the sum of the number of points
$t\in\overline{K}^\times$ where $E$ has multiplicative reduction and
twice the number of points $t\in\overline{K}^\times$ where $E$ has
additive reduction.  Clearly $0\le\deg\n-N^*(E)\le 4$ and so the
bound \ref{eq:geo-bound} gives $\rk_\Z E(F)\le N^*(E)$.

Now define $E_n$ as the elliptic curve defined by the equation of $E$
with $t$ replaced by $t^n$.  This is the base change of $E$ by the
the field homomorphism $K(t)\to K(t)$, $t\mapsto t^n$.  Clearly
$N^*(E_n)=nN^*(E)$ and so the geometric bound \ref{eq:geo-bound} gives
$\rk_\Z E_n(F)\le nN^*(E)$.

Assuming the Tate conjecture (namely the equality
$-\ord_{s=2}L(H^2(\EE),s)=\rk_\Z NS(\EE)$), Silverman proves by an
analytic method that
\begin{equation*}
\rk_\Z E_n(F)\le d_K(n)N^*(E)
\end{equation*}
where
\begin{equation*}
d_K(n)=\sum_{d|n}\frac{\phi(d)}{[K(\mu_d):K]}.
\end{equation*}
So, when $\mu_d\subset K$, $d_K(n)=n$ whereas if $K\cap\Q(\mu_n)=\Q$,
then $d_K(n)$ is the number of divisors of $n$.  Thus Silverman's
theorem gives an arithmetic bound for ranks of a very special class of
elliptic curves over function fields over number fields.  I believe
that there should be a much more general theorem in this direction.

There has been recent further work in this direction.  Namely,
Silverman \cite{Silverman2} has proven an interesting arithmetic bound
on ranks of elliptic curves over unramified, abelian towers, assuming
the Tate conjecture.  In the special case where the base curve is
itself elliptic and the tower is defined by the multiplication by $n$
isogenies, he obtains a very strong bound, stronger than what is
conjectured below.  (See his Theorem 2.)

Silverman also formulates a beautiful and precise conjecture along the
lines suggested above.  Namely, he conjectures that there is an
absolute constant $C$ such that for every non-isotrivial elliptic
curve over $F=K(\Curve)$ with conductor $\n$,
\begin{equation*}
\rk E(F)\mathrel{\mathop\le^?}C\frac{4g-4+\deg\n}{\log\deg\n}
\log|2\text{Disc}(K/\Q)|.
\end{equation*}
This conjecture is yet another instance of the fruitful interplay
between function fields and number fields.

\bibliographystyle{amsalpha}

\def\cprime{$'$} \def\cprime{$'$} \def\cprime{$'$} \def\cprime{$'$}
  \newcommand{\SortNoop}[1]{}
\providecommand{\bysame}{\leavevmode\hbox to3em{\hrulefill}\thinspace}
\providecommand{\MR}{\relax\ifhmode\unskip\space\fi MR }
\providecommand{\MRhref}[2]{%
  \href{http://www.ams.org/mathscinet-getitem?mr=#1}{#2}
}
\providecommand{\href}[2]{#2}

\end{document}